\documentclass[11pt]{amsart}

\usepackage[utf8]{inputenc} 	
\usepackage[T1]{fontenc}
\usepackage{lmodern}
\usepackage{yhmath}

\usepackage{amsthm,amsmath,amssymb,amsbsy,bbm,mathrsfs,supertabular,
eurosym,graphicx,enumitem,xcolor}
\usepackage{mathtools}

\newtheoremstyle{BBstyle0}  {}{}{\itshape}{}{\bfseries}{}{6pt}{}
\newtheoremstyle{BBstyle1}  {3pt}{3pt}{\rmfamily}{}{\itshape}{: }{3pt}{}
\newtheoremstyle{BBstyle2}  {3pt}{3pt}{\itshape}{}{\bfseries\large}{}{0pt}{}
\newtheoremstyle{BBstyle3}  {}{}{\itshape}{}{\bfseries}{: }{3pt}{}
\newtheoremstyle{BBstyle4}  {}{}{\rmfamily}{}{\bfseries}{}{6pt}{}
  
\usepackage[normalem]{ulem} 
\usepackage[authoryear]{natbib}

\newtheorem{thm}{Theorem}

\newtheorem{prop}{Proposition}
\newtheorem{df}{Definition}

\newtheorem{ass}{Assumption}

\theoremstyle{definition}

\usepackage[english]{babel}

\usepackage{enumitem}
\usepackage{hyperref}  
\newcommand{\pa}[1]{\left({#1}\right)}
\newcommand{\norm}[1]{\left\|{#1}\right\|}
\newcommand{\cro}[1]{\left[{#1}\right]}
\newcommand{\ab}[1]{\left|{#1}\right|}
\newcommand{\ac}[1]{\left\{{#1}\right\}}

\newcommand{\dfleche}[1]{\,\displaystyle{\mathop{\longrightarrow}_{#1}}\,}

\newcommand{\CV}[1]{\dfleche{#1}}

\newcommand{\E}{{\mathbb{E}}}

\renewcommand{\P}{{\mathbb{P}}}
 
\newcommand{\R}{{\mathbb{R}}}

\newcommand{\sF}{{\mathscr{F}}}

\DeclareMathAlphabet{\mathscrbf}{OMS}{mdugm}{b}{n}
\newcommand{\cE}{{\mathcal{E}}}
\newcommand{\cF}{{\mathcal{F}}}
\newcommand{\cG}{{\mathcal{G}}}

\newcommand{\cN}{{\mathcal{N}}}
\newcommand{\cO}{{\mathcal{O}}}

\newcommand{\cX}{{\mathcal{X}}}

\newcommand{\gf}{{\mathbf{f}}}
\renewcommand{\gg}{{\mathbf{g}}}

\newcommand{\gF}{{\mathbf{F}}}

\newcommand{\bs}[1]{\boldsymbol{#1}}

\newcommand{\bsY}{{\bs{Y}}} 
\newcommand{\bsZ}{{\bs{Z}}}
\newcommand{\frc}{{\mathfrak{c}}}

\newlist{lista}{enumerate}{1}
\setlist[lista,1]{label=\alph*),ref=\alph*)}

\newlist{listi}{enumerate}{1}
\setlist[listi,1]{label=(\roman*),ref=(\roman*),align=left}

\newcommand{\eref}[1]{(\ref{#1})}
\renewcommand{\ge}{\geqslant}
\renewcommand{\le}{\leqslant}
\newcommand{\1}{1\hskip-2.6pt{\rm l}}

\newcommand{\bigscal}[2]{\left\langle #1,#2\right\rangle}

\newcommand{\et}{^{\star}}

\newcommand{\ens}[1]{\{1,\ldots,#1\}}
\def\sgn{{\mathrm{\:sgn}}}

\DeclarePairedDelimiter\ceil{\lceil}{\rceil}
\DeclarePairedDelimiter\floor{\lfloor}{\rfloor}
\newcommand{\PES}[1]{\ceil*{#1}}
\newcommand{\PEI}[1]{\floor*{#1}}

\usepackage{pgfplots}
\pgfplotsset{compat=1.18}

\begin{document}
\title[]{Isotonic least squares is robust to heteroscedasticity, unimodal least squares is not.}
\author{Yannick BARAUD} 
\address{University of Luxembourg\\
DMATH\\
Maison du nombre\\
6 avenue de la Fonte\\
L-4364 Esch-sur-Alzette\\
Grand Duchy of Luxembourg.}
\email{yannick.baraud@uni.lu}

\date{\today}
\begin{abstract}
We compare the least squares estimator (LSE) with the $\ell_{1}$-estimator recently introduced by Baraud and Maillard~\citeyearpar{BM-2026} for estimating a regression function under a fixed design, when the errors are independent but possibly heteroscedastic and heavy-tailed. We show that the two estimators coincide on the class of nondecreasing (or nonincreasing) functions on an interval. As a consequence, the isotonic LSE inherits the non-asymptotic $\ell_{1}$-risk bounds of the $\ell_{1}$-estimator, which only require a moment of order $p\in[1,2]$ on the errors and depend on the average of their $p$-th absolute moments rather than on their maximum. This coincidence no longer holds on the nonconvex class of unimodal functions: we exhibit a heteroscedastic Gaussian example in which the LSE is inconsistent while the $\ell_{1}$-estimator converges at the parametric rate. Finally, we show on a one-dimensional linear model that neither estimator dominates the other: depending on the way the noise is distributed among the design points, the LSE may converge at rate $1/n$ while the $\ell_{1}$-estimator is inconsistent.
 \end{abstract}

\maketitle

\section{Introduction}
Let $x_{1,n},\ldots,x_{n,n}$ be $n$ deterministic points in a set $\cX$. For $i\in\ens{n}$, we associate with $x_{i,n}$ a response $Y_{i,n}\in\R$ and assume that these responses are (possibly) random, independent and integrable. From the observation of $\bsZ=\{(x_{1,n},Y_{1,n}),\ldots,(x_{n,n},Y_{n,n})\}$, we aim at estimating the function $f\et$ on $\cX$ that satisfies $f\et(x_{i,n})=\E\cro{Y_{i,n}}$ for $i\in\ens{n}$. The statistical model can also be written as 
\begin{equation}\label{eq-model}
Y_{i,n}=f\et(x_{i,n})+\sigma_{i,n}\xi_{i,n}\quad \text{for $i\in\{1,\ldots,n\}$}
\end{equation}
where the $\xi_{i,n}$ are independent centred random variables and the $\sigma_{i,n}$ are nonnegative numbers. This model includes various situations. In the literature,  the most classical one is that of homoscedastic regression{,} where the errors $\xi_{i,n}$ are independent (often i.i.d.) with variance 1 and the parameters $\sigma_{i,n}$ are constant. One may alternatively consider the case of independent Bernoulli or Poisson random variables $Y_{i,n}$ with their respective parameters $f\et(x_{i,n})$ for $i\in\ens{n}$. The statistical model then becomes heteroscedastic: the parameters $\sigma_{i,n}^{2}$ are given by $f\et(x_{i,n})(1-f\et(x_{i,n}))$ when the $Y_{i,n}$ are Bernoulli random variables and by $f\et(x_{i,n})$ when they are Poisson distributed.

To estimate $f\et$, we consider a class of functions $\cF$ on $\cX$, which we believe contains $f\et$ or at least provides a reasonable approximation of it, and we build an estimator $\widehat f$ with values in $\cF$.  

To solve this estimation problem, Baraud and Maillard~\citeyearpar{BM-2026} proposed an alternative to the classical Least Squares Estimator (LSE for short). We denote their estimator by $\widehat f$ and call it {\em an $\ell_{1}$-estimator} since they established non-asymptotic risk bounds for it with respect to the $\ell_{1}$-loss. These bounds do not require that the random variables $Y_{i,n}$ be square integrable but only that their centred $p$-absolute moments $\sigma_{i,n}^{p}=\sigma_{i,n}^{p}(p)=\E\cro{\ab{Y_{i,n}-\E\cro{Y_{i,n}}}^{p}}$ be finite for some $p\in [1,2]$. Furthermore, they do not depend on the maximum of the $\sigma_{i,n}^{p}(p)$ but on their average 
\begin{equation}\label{eq-sigmanp}
\sigma_{n}^{p}(p)=\frac{1}{n}\sum_{i=1}^{n}\E\cro{\ab{Y_{i,n}-\E\cro{Y_{i,n}}}^{p}}.
\end{equation}
This property ensures a certain stability of $\ell_{1}$-estimators with respect to a possible heteroscedasticity of the errors $Y_{i,n}-f\et(x_{i,n})$. Such a property is, in general, not satisfied {by} the LSE. Baraud and Maillard~\citeyearpar{BM-2026} provided an example in which $\cF$ is a linear space of dimension 1 containing $f\et$, the errors $Y_{i,n}-f\et(x_{i,n})$ are independent Gaussian random variables and the LSE is inconsistent for estimating $f\et$ while $\widehat f$ estimates it {at rate $\log(n)/\sqrt{n}$}. This example also shows that the LSE and the $\ell_{1}$-estimator do not necessarily coincide when $\cF$ is a linear space, hence, more generally, when $\cF$ is a convex set. 

Nevertheless, we show in this paper that they do coincide on the convex set of nondecreasing functions on an interval $\cX$. The same holds on the set of nonincreasing functions. As a consequence,  the LSE {inherits} the nice properties of the $\ell_{1}$-estimator on these classes. This is no longer the case when $\cF$ {now consists} of unimodal functions on an interval $\cX$. For such a class, we provide an example where the LSE becomes sensitive again to a possible heteroscedasticity of the errors. It may even be  inconsistent while the $\ell_{1}$-estimator estimates $f\et$ at {the} parametric rate. 

With the notable exception of  Efromovich and Pinsker~\citeyearpar{Efromovich1996} in the case where $\cF$ is a Sobolev space, the problem of estimating a regression function has almost exclusively been addressed through the study of the LSE. The case of shape-restricted regression has received special attention. 

In the case of isotonic regression, we refer the reader to Zhang~\citeyearpar{Zhang2002}, Chatterjee {\em et al{.}}~\citeyearpar{chatterjee2015}, Bellec~\citeyearpar{Bellec2018} {and} Guntuboyina and Sen~\citeyearpar{Guntuboyina2018}. 

The class $\sF$ of convex functions on an interval was considered in Chatterjee~\citeyearpar{Chatterjee2016} and Bellec~\citeyearpar{Bellec2018}. A more general study of the properties of the LSE on convex sets of functions can be found in Chatterjee~\citeyearpar{MR3269982}. 

For the study of the LSE on the (non-convex) set of unimodal functions, we refer to  Bellec~\citeyearpar{Bellec2018}, see also Chatterjee and Lafferty~\citeyearpar{Chatterjee2019}. Minami~\citeyearpar{Minami2020} considered the class of regression functions which are piecewise monotone under some constraints on the minimal length of the intervals which compose the partition. Feng {\em et al{.}}~\citeyearpar{S-shaped2022} studied the global rates and adaptation properties of the LSE for estimating a regression function with an $S$-shape on $[0,1]$, that is, a function which {is} first convex and then concave, the location of the inflexion point being unknown. Finally, let us mention the paper by Gao {\em et al{.}}~\citeyearpar{MR4102670} where the authors establish the minimax rates over classes of isotonic (and non-isotonic) piecewise constant signals with respect to the $\ell_{s}$-loss for $s\in [1,2]$. 

All these results assume that the errors are at least square integrable and {that} their variances {are} uniformly controlled (if not constant). 

The paper is organised as follows. In Section~\ref{sect-I}{,} we recall the construction of $\ell_{1}$-estimators and their properties. In Section~\ref{sect-II}, we establish that the LSE and the $\ell_{1}$-estimator coincide when $\cF$ consists of {nondecreasing} functions on an interval. The case of unimodal functions is tackled in Section~\ref{sect-III}. Finally, we show in Section~\ref{sect-IV} that when $\cF$ is a one{-}dimensional linear space, the LSE may perform much better than the $\ell_{1}$-estimator even in {a} heteroscedastic situation.

Throughout the paper, we use the following notation and conventions. The cardinality of a set $A$ is denoted by $|A|$ and we shall say that $A\subset \R$ is a {\em nontrivial interval} if its length is positive. Given a nonnegative number $x\in\R$, $\PES{x}$ is the {smallest integer which is at least as large as $x$} while $\PEI{x}$ is the largest integer not larger than $x$. We write $x\wedge y$ and $x\vee y$ respectively for the minimum and maximum {of} the two numbers $x$ and $y$ in $\R$ and
\[
\sgn(x)=\1_{x>0}-\1_{x<0}
=
\begin{cases}
1& \text{ if $x>0$},\\
0 &\text{ if $x=0$},\\
-1&\text{ if $x<0$}.\\
\end{cases}
\]
Finally, ${\bigscal{\cdot}{\cdot}}$ denotes the inner product of $\R^{n}$, $\ab{\cdot}$ the corresponding Euclidean norm and $\ab{\cdot}_{1}$ the $\ell_{1}$-norm.

\section{Construction of the estimator and its properties}\label{sect-I}
\subsection{Construction of the estimator}
The construction of the estimator is based on the test statistic $T(\bsZ,f,g)$ which is defined for  $f,g\in \cF$ by
\[
T(\bsZ,f,g)=\sum_{i=1}^{n}\pa{f(x_{i,n})-Y_{i,n}}\sgn\pa{f(x_{i,n})-g(x_{i,n})}\; \text{ and we set }\; T(\bsZ,f)=\sup_{g\in\cF}T(\bsZ,f,g).
\]
An $\ell_{1}$-estimator $\widehat f$ of $f\et$ {is defined as} an (approximate) minimiser of $f\mapsto T(\bsZ,f)$ over $\cF$, that is, an element of the (random) set
\begin{equation}\label{eq-defEk}
\cE_{\frc}(\bsZ,\cF)=\ac{f\in\cF,\; T(\bsZ,f)\le \inf_{g\in\cF}T(\bsZ,g)+\frc}
\end{equation}
where $\frc$ denotes some (small enough) nonnegative numerical constant. We shall say that $\widehat f$ is the $\ell_{1}$-estimator of $f\et$ when it is the {unique element of $\cE_{0}(\bsZ,\cF)$}.

More details about the construction and the properties of the test statistic $T$ are provided in Baraud and Maillard~\citeyearpar{BM-2026}.

\subsection{Extremal elements within a class $\cF$}
The elements $\overline f\in \cF$ which satisfy the following property play a central role in {Theorem~\ref{thm-main1} below}.
\begin{df}\label{def-extremal}
An element $\overline f$ in the class of functions $\cF$ on $\cX$ is extremal with degree not larger than $D_{n}\ge 1$ if the VC-dimensions of the two classes of subsets of $\cX_{n}=\{x_{1,n},\ldots,x_{n,n}\}$ given by 
\[
 \left\{\{f>\overline{f}\}\cap\cX_{n},\;  f\in \cF\right\}\quad \text{and}\quad \left\{ \{ f< \overline{f}\}\cap \cX_{n},\;  f \in \cF\right\}
\]
are not larger than $D_{n}$. The degree of an extremal element $\overline f\in{\cF}$ is the smallest
integer $D_{n}$ for which the above property is satisfied.
\end{df}
r{Since $\cX_{n}$ has at most $n$ elements, every element of $\cF$ is extremal with degree not larger than $n$.} It is shown in Baraud and Maillard~\citeyearpar{BM-2026} that if $\cF$ is a linear space of dimension $d$, all the elements of $\cF$ are extremal with degree not larger than $d+1$. If $\cX$ is a nontrivial interval and $\cF$ {is} the class of monotone functions on $\cX$, the elements of $\cF$ which are piecewise constant on a partition of $\cX$ into at most $d$ intervals are extremal with degree $D_{n}\le d+1$. If, on this interval, $\cF$ consists of functions which are convex or concave, then the elements of $\cF$ which are piecewise linear and based on a partition of $\cX$ into at most $d$ intervals are extremal with degree not larger than $5+2d$. More examples can be found in Baraud and Maillard~\citeyearpar{BM-2026}. Hereafter, for every $D\in\{1,\ldots,n\}$, we denote by $\cO(D)$ the subset of $\cF$ that consists of the functions which are extremal with a degree not larger than $D$.

\subsection{The $\ell_{1}$-risk of the estimator}
Given two functions $f,g$ on $\cX$ and a subset $\cG$ of $\cF$, we set for $s\in\{1,2\}$
\[
\ell_{s}(f,g)=\frac{1}{n}\sum_{i=1}^{n}\ab{f(x_{i,n})-g(x_{i,n})}^{s}\quad \text{and}\quad \ell_{s}(f,\cG)=\inf_{g\in\cG}\ell_{s}(f,g).
\]

The following result is proven in  Baraud and Maillard~\citeyearpar{BM-2026}.
\begin{thm}\label{thm-main1}
Assume that the random variables $Y_{1,n},\ldots,Y_{n,n}$ are independent and let $\widehat f$ be an element of $\cE_{\frc}(\bsZ,\cF)$ with $\frc\ge 0$. Then, for every $\overline f\in\cO(D)$ with $D\in\{1,\ldots,n\}$,
\begin{equation}\label{eq-main00}
\E\cro{\ell_{1}\pa{\widehat f,\overline f}}\le 2\ell_{1}\pa{f\et,\overline f}+ \kappa\inf_{p\in [1,2]}\cro{\sigma_{n}(p)\pa{\frac{D}{n}}^{1-\frac{1}{p}}}+\frac{\frc}{n}
\end{equation}
where $\kappa$ denotes a positive numerical constant and $\sigma_{n}(p)$ is defined by \eref{eq-sigmanp}. In particular,
\begin{align}
\E\cro{\ell_{1}\pa{f\et,\widehat f}}\le  \inf_{D\in \{1,\ldots,n\}}\ac{3\ell_{1}\pa{f\et,\cO(D)}+ \kappa\inf_{p\in [1,2]}\cro{\sigma_{n}(p)\pa{\frac{D}{n}}^{1-\frac{1}{p}}}}+\frac{\frc}{n}.\label{eq-main01}
\end{align}
%
\end{thm}
As a consequence of \eref{eq-main00}, when the $\sigma_{i,n}=\sigma>0$ are constant and the $\xi_{i,n}$ are i.i.d.\ random variables with variance 1, an $\ell_{1}$-estimator {converges at the parametric rate $1/\sqrt{n}$ when $f\et$ is an extremal element of $\cF$} whose degree is bounded independently of $n$. For $p\in (1,2)$, this rate becomes $n^{-(1-1/p)}$ and  Devroye {\em et al.}~\citeyearpar{MR3576558} showed that it cannot be improved in general when the $\xi_{i,n}$ are i.i.d.\ with only a moment of order $p$.

Inequality \eref{eq-main01} shows a certain stability of the risk of an $\ell_{1}$-estimator. If the regression function $f\et$ (which might not belong to $\cF$) is close enough to an extremal function $\overline f\in\cF$, the risk bound that we get remains of the same order as the one we would get if the regression function were $\overline f$.

\section{The main result}\label{sect-II}
Let $x_{1,n}\le \ldots\le x_{n,n}$ be $n$ points of a nontrivial interval $\cX$ and $\cF_{\uparrow}$ the class of nondecreasing functions on $\cX$.%
\begin{thm}\label{thm-mono}
For $\cF=\cF_{\uparrow}$, the LSE $\widetilde f$ minimises $f\mapsto T(\bsZ,f)$ over $\cF$ and even satisfies $T(\bsZ,\widetilde f)=0$. In particular, the LSE satisfies \eref{eq-main00} and \eref{eq-main01} with $\frc=0$. If $\widehat f$ is another minimiser of $f\mapsto T(\bsZ,f)$ then $\widehat f=\widetilde f$ on $\{x_{1,n},\ldots,x_{n,n}\}$.
\end{thm}

In our framework, it is convenient to identify a function $f$ on $\cX$ with the vector $\gf=(f(x_{1,n}),\ldots,f(x_{n,n}))$ in $\R^{n}$, using a bold letter for the corresponding vector. Then, a class of functions $\cF$ is associated with the set of vectors $\gF=\{\gf, f \in\cF\}\subset \R^{n}$ and we may write our criterion as 
\[
T(\bsZ,\gf)=\sup_{\gg\in\gF}T(\bsZ,\gf,\gg)\quad \text{with}\quad T(\bsZ,\gf,\gg)=\bigscal{\gf-\bsY}{\sgn(\gf-\gg)}\quad \text{for $\gf\in\gF$,}
\]
where $\bsY=(Y_{1,n},\ldots,Y_{n,n})$ and $\sgn(\gf-\gg)=(\sgn(f(x_{1,n})-g(x_{1,n})),\ldots,\sgn(f(x_{n,n})-g(x_{n,n})))$. When $\cF=\cF_{\uparrow}$, {note} that the LSE $\widetilde f$ is associated with the projection $\widetilde \gf$ of $\bsY$ onto the closed convex set $\gF_{\uparrow}\subset \R^{n}$. {Theorem~\ref{thm-mono} states} that $\widetilde \gf$ is the unique minimiser {of} $\gf\mapsto T(\bsZ,\gf)$ on $\gF_{\uparrow}$. 

The conclusion of Theorem~\ref{thm-mono} extends to {any class} of functions $\cF$ {satisfying the following assumption}:
\begin{ass}\label{Ass1}
The set $\gF$ is a closed convex subset of $\R^{n}$ and for every $\gf, \gg\in\gF$, there exists $\lambda>0$ such that $\gf-\lambda\sgn(\gf-\gg)$ belongs to $\gF$.
\end{ass}

\begin{proof}[Proof of Theorem~\ref{thm-mono}]
As already mentioned, $\widetilde \gf$ is {the} projection of $\bsY$ onto the closed convex set $\gF=\gF_{\uparrow}$. It therefore exists and is unique. It is characterised by the fact that
\begin{equation}\label{eq-LS}
\bigscal{\bsY-\widetilde \gf}{\gg-\widetilde \gf}\le 0\quad \text{for all $\gg\in\gF$.}
\end{equation}
Note that for every $\gf\in\gF$, $T(\bsZ,\gf)\ge T(\bsZ,\gf,\gf)=\bigscal{\gf-\bsY}{\sgn(\gf-\gf)}=0$. It is therefore sufficient to show that $T(\bsZ,\widetilde \gf)=0$ in order to prove that $\widetilde \gf$ minimises $\gf\mapsto T(\bsZ,\gf)$ on $\gF$. The result is clear when $n=1$ since ${\widetilde \gf}=\bsY=Y_{1,1}$. We may therefore assume $n\ge 2$ hereafter. By definition of $T(\bsZ,{\widetilde \gf})$, it is sufficient to prove that for every $\gg\in\gF$, $T({\bsZ},\widetilde \gf,\gg)\le 0$. Provided that Assumption~\ref{Ass1} is satisfied, we may choose in \eref{eq-LS} $\widetilde \gg^{\lambda}=\widetilde \gf-\lambda\sgn(\widetilde \gf-\gg)$ in place of $\gg$ and obtain that 
\[
0\ge \bigscal{\bsY-\widetilde \gf}{\widetilde \gg^{\lambda}-\widetilde \gf}=\lambda\bigscal{\widetilde \gf-\bsY}{\sgn(\widetilde \gf-\gg)}=\lambda T(\bsZ,\widetilde \gf,\gg)
\]
which gives $T(\bsZ,\widetilde \gf,\gg)\le 0$ since $\lambda>0$. It remains thus to check that Assumption~\ref{Ass1} is satisfied. 
Given $\gf,\gg\in\gF$ and $\lambda>0$, we set $\gg^{\lambda}=\gf-\lambda\sgn(\gf-\gg)$ and note that for every $i\in\{2,\ldots,n\}$, 
\begin{align*}
g^{\lambda}(x_{i,n})-g^{\lambda}(x_{i-1,n})={f(x_{i,n})-f(x_{i-1,n})}+\lambda \Delta_{i}
\end{align*}
with 
\[
\Delta_{i}=\sgn(f(x_{i-1,n})-g(x_{i-1,n}))-\sgn(f(x_{i,n})-g(x_{i,n})).
\]
Let us now distinguish between two cases. If $f(x_{i,n})>f(x_{i-1,n})$, $g^{\lambda}(x_{i,n})-g^{\lambda}(x_{i-1,n})\ge 0$ as soon as $0<\lambda\le (f(x_{i,n})-f(x_{i-1,n}))/2$ since $|\Delta_{i}|\le 2$. Otherwise, $f(x_{i,n})=f(x_{i-1,n})$ and since $g(x_{i-1,n})\le g(x_{i,n})$, $\Delta_{i}$ is nonnegative and so is $g^{\lambda}(x_{i,n})-g^{\lambda}(x_{i-1,n})$ for every value of $\lambda>0$. We conclude that for the choice 
\[
\lambda=\min\ac{\frac{f(x_{i,n})-f(x_{i-1,n})}{2},\; i\in\{2,\ldots,n\} \text{ with } f(x_{i,n})>f(x_{i-1,n})},
\]
and $\lambda=1$ when $f$ is constant on $\cX_{n}$,  the vector $\gg^{\lambda}$ belongs to $\gF$ and Assumption~\ref{Ass1} is satisfied. 

Let us now show that $\widetilde \gf$ is the only minimiser of $\gf\mapsto T({\bsZ},\gf)$. Let $\widehat \gf$ in $\gF$ be such that $T({\bsZ},\widehat \gf)=T({\bsZ},\widetilde \gf)=0$. This implies that $T({\bsZ},\widehat \gf,\widetilde \gf)\le 0$ and $T({\bsZ},\widetilde \gf,\widehat \gf)\le 0$. 
In particular, 
\begin{align*}
0&\ge T({\bsZ},\widehat \gf,\widetilde \gf)+T({\bsZ},\widetilde \gf,\widehat \gf)=\bigscal{\widehat \gf-\bsY}{\sgn\pa{\widehat \gf-\widetilde \gf}}+\bigscal{\widetilde \gf-\bsY}{\sgn\pa{\widetilde \gf-\widehat \gf}}\\
&={\bigscal{\widehat \gf-\widetilde \gf}{\sgn\pa{\widehat \gf-\widetilde \gf}}}=|\widetilde \gf-\widehat \gf|_{1}{.}
\end{align*}
Hence $\widehat \gf=\widetilde \gf$.
\end{proof}

\section{The case of unimodal functions on an interval}\label{sect-III}
In {this} section{,} we consider the class $\cF$ of functions $f$ which are unimodal on the nontrivial interval $\cX$. This means that either $f$ is monotone on $\cX$, or there exists $m$ in the interior of $\cX$ such that $f$ is nondecreasing on $(-\infty, m]\cap \cX$ and nonincreasing on $[m,+\infty)\cap \cX$. In this case, we shall see that the LSE and the $\ell_{1}$-estimator may not coincide. 

For illustration, let us consider Model~\eref{eq-model} where $n\ge 8$, the $\xi_{i,n}$ are {i.i.d.\ standard} Gaussian random variables, $x_{i,n}=i/n\in\cX=[0,1]$ for $i\in\ens{n}$, $f\et=\1_{[1/4,1/2]}\in\cF$ and $\sigma_{i,n}=0$ for all $i\in\ens{n}$ except for $i=i_{0}=\PES{3n/4}$ for which we set $\sigma_{i_{0},n}=\sqrt{n}$. Then the following holds.
\begin{thm}\label{thm-unimodal}
Let $\widetilde f_{n}$ and $\widehat f_{n}$ be respectively a LSE and an $\ell_{1}$-estimator of $f\et$ on the class of unimodal functions on $[0,1]$ in the statistical model described above. There exists an {event} $\Omega_{n}$, whose probability is at least 30\% {for $n$ large enough}, on which 
\[
\sqrt{\ell_{2}(\widetilde f_{n},f\et)}\ge \ell_{1}(\widetilde f_{n},f\et)\ge \frac{1}{8}\pa{1-\frac{4}{n}}.
\]
In particular, $\widetilde f_{n}$ is not consistent while $\widehat f_{n}$ satisfies
\begin{equation}\label{eq-fhat}
\E\cro{\ell_{1}(\widehat f_{n},f\et)}\le \frac{\kappa\sqrt{6}}{\sqrt{n}}+\frac{\frc}{n}.
\end{equation}
\end{thm}
Inequality \eref{eq-fhat} is a consequence of \eref{eq-main00}. In the statistical model we consider here, $f\et$ is piecewise constant on a partition of $\cX$ into $K=3$ pieces and is therefore extremal in $\cF$ with degree not larger than $D_{n}=2k+K-1=6$, by applying Proposition~4.1 in Baraud and Maillard~\citeyearpar{BM-2026} with $k=2$ and $K=3$. Besides, $\sigma_{n}(2)=1$ and the result follows by applying Theorem~\ref{thm-main1} with $\overline f=f\et\in\cO(6)$.

\begin{proof}[Proof of Theorem~\ref{thm-unimodal}]
Note that $Y_{i,n}=f\et(x_{i,n})$ for $i\ne i_{0}$. We define $E_{n}=\{1,\ldots,n\}$, 
$I_{1}=\{i\in E_{n},\ x_{i,n}\in [0,1/4)\}$, $I_{2}=\{i\in E_{n},\ x_{i,n}\in [1/4,1/2]\}$, $I_{3}=\{i\in E_{n},\ x_{i,n}\in (1/2,i_{0}/n)\}$, $I_{4}=\{i\in E_{n},\ x_{i,n}\in (i_{0}/n,1]\}$. For $n\ge 8$, one can check that
\[
|I_{j}|\ge N=\frac{n}{4}-1\ge 1\quad \text{for every $j\in \{1,\ldots,4\}$ and }\quad |I_{3}|\le \frac{n}{4}+1.
\]
We denote by $\overline f$ the function which satisfies $\overline f(x_{i,n})=|I_{2}|/[|I_{2}|+|I_{3}|]\in (0,1)$ for $i\in I_{2}\cup I_{3}$, $\overline f(x_{i_{0},n})=Y_{i_{0},n}$, $\overline f(x_{i,n})=0$ for $i\in I_{1}\cup I_{4}$ and $\overline f$ is affine between {two} consecutive points $x_{i,n}$ and takes the value 0 elsewhere. The function $\overline f$ belongs to $\cF$ on the event $\Omega_{n}=\{Y_{i_{0},n}\ge 1/2+\sqrt{n/4+1}\}$ and we note that 
\begin{align*}
\sum_{i=1}^{n}\pa{Y_{i,n}-\overline f (x_{i,n})}^{2}&=\sum_{i\in I_{2}}\pa{1-\frac{|I_{2}|}{|I_{2}|+|I_{3}|}}^{2}+\sum_{i\in I_{3}}\pa{\frac{|I_{2}|}{|I_{2}|+|I_{3}|}}^{2}=\frac{|I_{2}||I_{3}|}{|I_{2}|+|I_{3}|}\\
&\le |I_{3}|\le \frac{n}{4}+1.
\end{align*}
Let $z=1/2$ and $f$ be an element of $\cF$ that satisfies $\sum_{i\ne i_{0}}\ab{f(x_{i,n})-f\et(x_{i,n})}< zN$. Then for every $j\in \{1,\ldots,4\}$, $\sum_{i\in I_{j}}\ab{f(x_{i,n})-f\et(x_{i,n})}< zN$ and since $|I_{j}|\ge N$, there exists $i_{j}\in I_{j}$ such that $\ab{f(x_{i,n})-f\et(x_{i,n})}< z$. Hence, $f(x_{i_{1},n})< z$, $f(x_{i_{2},n})> 1-z$, $f(x_{i_{3},n})< z$ and  $f(x_{i_{4},n})< z$. The mode of $f$ is therefore located in the interval $(x_{i_{1},n},x_{i_{3},n})$. Since $x_{i_{0},n}>x_{i_{3},n}$, we necessarily have $f(x_{i_{0},n})<z$. We deduce that if $f\in\cF$ satisfies $\sum_{i\ne i_{0}}\ab{f(x_{i,n})-f\et(x_{i,n})}< zN$, {then,} on the {event} $\Omega_{n}${,}
\[
\sum_{i=1}^{n}\pa{Y_{i,n}-f(x_{i,n})}^{2}\ge \pa{Y_{i_{0},n}-f(x_{i_{0},n})}^{2}>  \pa{Y_{i_{0},n}-\frac{1}{2}}^{2}\ge |I_{3}|.
\]
Hence, the function $f$ cannot be the LSE{. Consequently, $\sum_{i\ne i_{0}}\ab{\widetilde f_{n}(x_{i,n})-f\et(x_{i,n})}\ge zN$ and} $\widetilde f_{n}$ necessarily satisfies on $\Omega_{n}$
\begin{align*}
\ell_{2}\pa{\widetilde f_{n},f\et}\ge \pa{\frac{1}{n}\sum_{i=1}^{n}\ab{\widetilde f_{n}(x_{i,n})-f\et(x_{i,n})}}^{2}\ge \frac{z^{2}N^{2}}{n^{2}}=\frac{1}{4}\pa{\frac{1}{4}-\frac{1}{n}}^{2}.
\end{align*}
In particular{,} $\ell_{1}(\widetilde f_{n},f\et)\ge (1-4/n)/8$ on $\Omega_{n}$ and since $Y_{i_{0},n}$ has the same distribution as $\sqrt{n}Z$ with $Z\sim\cN(0,1)$, we obtain that 
\[
\lim_{n\to +\infty}\P\pa{\Omega_{n}}=\lim_{n\to +\infty}\P\cro{Z\ge \frac{1}{2\sqrt{n}}+\sqrt{\frac{1}{4}+\frac{1}{n}}}=\P\cro{Z\ge \frac{1}{2}}\ge 30\%.
\]
\end{proof}

\section{The LSE and the $\ell_{1}$-estimator on a one-dimensional linear space}\label{sect-IV}
The previous comparison of the LSE with the $\ell_{1}$-estimator should not lead to the conclusion that the latter performs better than the former in heteroscedastic regression models. In fact, there are examples where the LSE {converges} (quite fast) to $f\et$ although the $\ell_{1}$-estimator is not even consistent. This is what we illustrate here when $\cF$ is the linear span of a function $f_{1}$ on $\cX=(0,1]$. In this section, we consider Model~{\eref{eq-model}} where $f\et=a\et f_{1}$, $x_{i,n}=i/n$ for $i\in\{1,\ldots,n\}$ and {the} $\xi_{i,n}$ are i.i.d.\ standard Gaussian random variables. We assume that $f_{1}$ is positive and continuous on $\cX$ and satisfies
\[
\norm{f_{1}}_{1,n}=\frac{1}{n}\sum_{i=1}^{n}f_{1}(x_{i,n})\CV{n\to +\infty}\int_{0}^{1}f_{1}(x)dx
\]
and 
\[
\norm{f_{1}}_{2,n}^{2}=\frac{1}{n}\sum_{i=1}^{n}f_{1}^{2}(x_{i,n})\CV{n\to +\infty}\int_{0}^{1}f_{1}^{2}(x)dx<+\infty.
\]
For every $a\in\R$, we set $f_{a}=af_{1}\in\cF$ so that 
\[
\ell_{1}(f_{a},f_{a\et})=\ab{a-a\et}\norm{f_{1}}_{1,n}\quad \text{while}\quad \ell_{2}(f_{a},f_{a\et})=\ab{a-a\et}^{2}\norm{f_{1}}_{2,n}^{2}.
\]
In particular an estimator $f_{\overline a_{n}}$ of $f\et=f_{a\et}$ is consistent in the $\ell_{s}$-loss with $s\in\{1,2\}$ if and only if $\overline a_{n}$ is a consistent estimator of $a\et$. 

\begin{prop}\label{prp-comp}
In the statistical model described above, the LSE $f_{\widetilde a_{n}}$ and the $\ell_{1}$-estimator $f_{\widehat a_{n}}$ exist and are unique. They are defined by the values 
\[
\widetilde a_{n}=\frac{\sum_{i=1}^{n}Y_{i,n}f_{1}(x_{i,n})}{n{\norm{f_{1}}_{2,n}^{2}}}\quad \text{and}\quad \widehat a_{n}=\frac{\sum_{i=1}^{n}Y_{i,n}}{n\norm{f_{1}}_{1,n}}.
\]
Besides,
\begin{equation}\label{eq-comp}
 \norm{f_{1}}_{2,n}^{4}\E\cro{\ab{\widetilde a_{n}-a\et}^{2}}=\frac{\sum_{i=1}^{n}\sigma_{i,n}^{2}f_{1}^{2}(x_{i,n})}{n^{2}},\quad \norm{f_{1}}_{1,n}^{2}\E\cro{\ab{\widehat a_{n}-a\et}^{2}}=\frac{\sum_{i=1}^{n}\sigma_{i,n}^{2}}{n^{2}}.
\end{equation}
In particular,
\[
\min_{i\in\{1,\ldots,n\}}f_{1}^{2}(x_{i,n})\le \frac{\norm{f_{1}}_{2,n}^{4}}{\norm{f_{1}}_{1,n}^{2}}\frac{\E\cro{\ab{\widetilde a_{n}-a\et}^{2}}}{\E\cro{\ab{\widehat a_{n}-a\et}^{2}}}= \frac{\sum_{i=1}^{n}\sigma_{i,n}^{2}f_{1}^{2}(x_{i,n})}{\sum_{i=1}^{n}\sigma_{i,n}^{2}}\le \max_{i\in\{1,\ldots,n\}}f_{1}^{2}(x_{i,n}).
\]
\end{prop}
We deduce from this proposition that when the function $f_{1}$ is bounded away from 0 and infinity on $\cX=(0,1]${,} the quadratic risks of $\widetilde a_{n}$ and $\widehat a_{n}$ are of the same order. This may no longer be the case when $f_{1}$ can be close to 0 or tend to infinity in the neighbourhood of 0. For example, if $f_{1}(x)=x$ on $(0,1]$, $\sigma_{1,n}=\tau_{n}$ and $\sigma_{i,n}=0$ for every $i\ne 1$, then $\min_{i\in\{1,\ldots,n\}}f_{1}^{2}(x_{i,n})=n^{-2}$, 
\[
\E\cro{\ab{\widetilde a_{n}-a\et}^{2}}=\frac{\norm{f_{1}}_{1,n}^{2}}{\norm{f_{1}}_{2,n}^{4}}\frac{\E\cro{\ab{\widehat a_{n}-a\et}^{2}}}{n^{2}}=\frac{\tau_{n}^{2}}{n^{4}\norm{f_{1}}_{2,n}^{4}}.
\]
Not only {does} the quadratic risk of $\widetilde a_{n}$ {improve on} that of $\widehat a_{n}$ by a factor of order $1/n^{2}${,} but $\widehat a_{n}$ may even be inconsistent for $\tau_{n}=n$ while $\widetilde a_{n}$ still converges at rate $1/n$. 

The opposite situation may also occur if now $f_{1}$ becomes large in the neighbourhood of 0 as in the example provided in {Section~2.3 of Baraud and Maillard~\citeyearpar{BM-2026}}.

\begin{proof}
Let us set $f_{a}=af_{1}\in\cF$ for $a\in\R$. Since $f_{1}$ is positive on $(0,1]$, for every $f_{a},f_{b}\in {\cF}$ with $a,b\in\R$, 
\begin{align*}
T(\bsZ,f_{a},f_{b})&=\sum_{i=1}^{n}\pa{f_{a}(x_{i,n})-Y_{i,n}}\sgn\pa{af_{1}(x_{i,n})-bf_{1}(x_{i,n})}\\
&=\sgn(a-b)\sum_{i=1}^{n}\pa{af_{1}(x_{i,n})-Y_{i,n}}.
\end{align*}
Hence
\[
T(\bsZ,f_{a})=\sup_{b\in\R}T(\bsZ,f_{a},f_{b})=\ab{\sum_{i=1}^{n}\pa{af_{1}(x_{i,n})-Y_{i,n}}}
\]
and the minimum is achieved for 
\[
\widehat a_{n}=\frac{\sum_{i=1}^{n}Y_{i,n}}{\sum_{i=1}^{n}f_{1}(x_{i,n})}=a\et+\frac{\sum_{i=1}^{n}\sigma_{i,n}\xi_{i,n}}{\sum_{i=1}^{n}f_{1}(x_{i,n})}=a\et+\frac{\sum_{i=1}^{n}\sigma_{i,n}\xi_{i,n}}{n\norm{f_{1}}_{1,n}}.
\]
The function $f_{\widehat a_{n}}$ is therefore the $\ell_{1}$-estimator while the LSE is $f_{\widetilde a_{n}}$ with 
\[
\widetilde a_{n}=\frac{\sum_{i=1}^{n}Y_{i,n}f_{1}(x_{i,n})}{\sum_{i=1}^{n}f^{2}_{1}(x_{i,n})}=a\et+\frac{\sum_{i=1}^{n}\sigma_{i,n}f_{1}(x_{i,n})\xi_{i,n}}{\sum_{i=1}^{n}f^{2}_{1}(x_{i,n})}=a\et+\frac{\sum_{i=1}^{n}\sigma_{i,n}f_{1}(x_{i,n})\xi_{i,n}}{n\norm{f_{1}}_{2,n}^{2}}
\]
and the equalities in \eref{eq-comp} follow.

\end{proof}


\begin{thebibliography}{}

\bibitem[Baraud and Maillard, 2026]{BM-2026}
Baraud, Y. and Maillard, G. (2026).
\newblock Estimating a regression function under possible heteroscedastic and
  heavy-tailed errors. {Application} to shape-restricted regression.
\newblock {\em to appear in Math. Stat. Learn.}

\bibitem[Bellec, 2018]{Bellec2018}
Bellec, P.~C. (2018).
\newblock {Sharp oracle inequalities for Least Squares estimators in shape
  restricted regression}.
\newblock {\em The Annals of Statistics}, 46(2):745 -- 780.

\bibitem[Chatterjee, 2014]{MR3269982}
Chatterjee, S. (2014).
\newblock A new perspective on least squares under convex constraint.
\newblock {\em Ann. Statist.}, 42(6):2340--2381.

\bibitem[Chatterjee, 2016]{Chatterjee2016}
Chatterjee, S. (2016).
\newblock {An improved global risk bound in concave regression}.
\newblock {\em Electronic Journal of Statistics}, 10(1):1608 -- 1629.

\bibitem[Chatterjee et~al., 2015]{chatterjee2015}
Chatterjee, S., Guntuboyina, A., and Sen, B. (2015).
\newblock On risk bounds in isotonic and other shape restricted regression
  problems.
\newblock {\em Ann. Statist.}, 43(4):1774--1800.

\bibitem[Chatterjee and Lafferty, 2019]{Chatterjee2019}
Chatterjee, S. and Lafferty, J. (2019).
\newblock {Adaptive risk bounds in unimodal regression}.
\newblock {\em Bernoulli}, 25(1):1 -- 25.

\bibitem[Devroye et~al., 2016]{MR3576558}
Devroye, L., Lerasle, M., Lugosi, G., and Oliveira, R.~I. (2016).
\newblock Sub-{G}aussian mean estimators.
\newblock {\em Ann. Statist.}, 44(6):2695--2725.

\bibitem[Efromovich and Pinsker, 1996]{Efromovich1996}
Efromovich, S. and Pinsker, M. (1996).
\newblock Sharp-optimal and adaptive estimation for heteroscedastic
  nonparametric regression.
\newblock {\em Statistica Sinica}, 6:925--942.

\bibitem[Feng et~al., 2022]{S-shaped2022}
Feng, O.~Y., Chen, Y., Han, Q., Carroll, R.~J., and Samworth, R.~J. (2022).
\newblock Nonparametric, tuning-free estimation of {{S}}-shaped functions.
\newblock {\em Journal of the Royal Statistical Society Series B: Statistical
  Methodology}, 84(4):1324--1352.

\bibitem[Gao et~al., 2020]{MR4102670}
Gao, C., Han, F., and Zhang, C.-H. (2020).
\newblock On estimation of isotonic piecewise constant signals.
\newblock {\em Ann. Statist.}, 48(2):629--654.

\bibitem[Guntuboyina and Sen, 2018]{Guntuboyina2018}
Guntuboyina, A. and Sen, B. (2018).
\newblock {Nonparametric Shape-Restricted Regression}.
\newblock {\em Statistical Science}, 33(4):568 -- 594.

\bibitem[Minami, 2020]{Minami2020}
Minami, K. (2020).
\newblock {Estimating piecewise monotone signals}.
\newblock {\em Electronic Journal of Statistics}, 14(1):1508 -- 1576.

\bibitem[Zhang, 2002]{Zhang2002}
Zhang, C.-H. (2002).
\newblock {Risk bounds in isotonic regression}.
\newblock {\em The Annals of Statistics}, 30(2):528 -- 555.

\end{thebibliography}

\end{document}